\documentclass[11pt]{article}

\usepackage{amsmath,amsfonts,amssymb,amsthm}
\usepackage[margin=1in]{geometry}
\usepackage[colorlinks]{hyperref}
\usepackage{enumitem}
\usepackage{mathtools}
\usepackage{enumitem}
\usepackage{algorithm, caption}
\usepackage[noend]{algpseudocode}
\usepackage{subfigure}
\usepackage{graphicx}
\usepackage{tabularx,environ}
\usepackage[braket]{qcircuit}
\usepackage[dvipsnames]{xcolor}
\usepackage{multirow}
\usepackage[capitalize]{cleveref}
\usepackage{dsfont}
\usepackage{comment}

\makeatletter
\newcommand*\bigcdot{\mathpalette\bigcdot@{.5}}
\newcommand*\bigcdot@[2]{\mathbin{\vcenter{\hbox{\scalebox{#2}{$\m@th#1\bullet$}}}}}
\makeatother

    \hypersetup{
    	colorlinks,
    	linkcolor={red!100!black},
    	citecolor={blue!100!black},
    }
\makeatletter
\newcommand\notsotiny{\@setfontsize\notsotiny\@vipt\@viipt}
\makeatother

\newcommand{\norm}[1]{\left\|{#1}\right\|}

\def\01{\{0,1\}}

\theoremstyle{plain}

\newtheorem{theorem}{Theorem}

\newtheorem{fact}[theorem]{Fact}


\definecolor{applegreen}{rgb}{0.0, 0.5, 0.0}

\newcommand{\Id}{\ensuremath{\mathop{\rm Id}\nolimits}}

\DeclareMathOperator{\Tr}{Tr}

\DeclareMathOperator{\rank}{rk}

\newcommand{\beq}{\begin{equation}}
\newcommand{\eeq}{\end{equation}}
\newcommand{\beqn}{\begin{equation*}}
\newcommand{\eeqn}{\end{equation*}}
\newcommand{\beqr}{\begin{eqnarray}}
\newcommand{\eeqr}{\end{eqnarray}}
\newcommand{\beqrn}{\begin{eqnarray*}}
\newcommand{\eeqrn}{\end{eqnarray*}}
\newcommand{\bmline}{\begin{multline}}
\newcommand{\emline}{\end{multline}}
\newcommand{\bmlinen}{\begin{multline*}}
\newcommand{\emlinen}{\end{multline*}}

\theoremstyle{plain}

\newtheorem{question}[theorem]{Question}

\theoremstyle{definition}

\theoremstyle{remark}

\renewenvironment{proof}[1][]{
    \begin{trivlist}
    \item[\hspace{\labelsep}{\em\noindent Proof#1:\/}]}
    {{\hfill$\Box$}
    \end{trivlist}
}

\newtheoremstyle{named}{}{}{\itshape}{}{\bfseries}{.}{.5em}{\thmnote{#3}}
\theoremstyle{named}

\usepackage{tikz}
\usepackage{thm-restate} %To repeat theorems, lemmas etc....
\title{All $S_p$ notions of quantum expansion are equivalent}
\author{ 
Francisco Escudero Gutiérrez\thanks{Qusoft and CWI \href{feg@cwi.nl}{feg@cwi.nl} } \and Garazi Muguruza\thanks{Qusoft and UvA \href{g.muguruzalasa@uva.nl}{g.muguruzalasa@uva.nl} } }

\date{  }

\begin{document}

\maketitle
\begin{abstract}
    In a recent work Li, Qiao, Wigderson, Wigderson and Zhang introduced notions of quantum expansion based on $S_p$ norms and posed as an open question if they were all equivalent. We give an affirmative answer to this question. 
\end{abstract}

\section{Introduction}
In this note we consider the problem of showing that certain notions of quantum expansion are equivalent. Quantum expansion is an analogue of the well-studied concept of graph expansion. The latter has many equivalent versions that let it play an important role in many fields such as geometry, group-theory, combinatorics and theoretical computer science (for a survey see \cite{hoory2006expander}). 

Given a $d$-regular graph $G=([n],E)$, its  \emph{edge expansion} is given by $$h(G)=\min_{\substack{W\subset [n],|W|\leq n/2}} \frac{|\partial W|}{d|W|},$$
where $|W|$ stands for the size of $W$ and $|\partial W|$ is the number of edges that go from $W$ to its complement. Let  $\mathbf{B}=(B_1,\dots, B_d)$ be permutation matrices  that decompose the adjacency matrix $A$ of $G$ as $\sum_{i\in [n]}{B_i}$\footnote{It follows from Hall's theorem that this decomposition always exists.}. In particular, $\mathbf B$ is a bistochastic tuple, meaning that $\sum_{i\in[n]}B_i^*B_i=\sum_{i\in[n]}B_iB_i^*=d\cdot\Id$. Thus, it induces a quantum unital  channel via
$$\Phi_{\mathbf B}(X)=\frac{1}{d}\sum_{i\in [n]}B_iXB_i^*.$$
A quantum channel is a linear operator from $M_n(\mathbb C)$ to $M_n(\mathbb C)$ that maps quantum states (positive semidefinite matrices with trace 1) to quantum states. It is unital if it maps the identity to itself. 
Note that if $\mathbf B$ is a permutation tuple, $W$ a subset of $[n]$ and $P_W$ the projector onto the space spanned by the canonical basis vectors indicated by $W$, then $\langle \Id-P_W,\Phi_{\mathbf B}(P_W)\rangle=|\partial W|$. Here $\langle \cdot,\cdot\rangle$ stands for the inner product $\langle M,N\rangle=\Tr[M^*N]$.  Motivated by this, Hastings introduced the \emph{quantum edge expansion} of a bistochastic tuple $\mathbf B\in M_n(\mathbb C)^d$, which is given by 
\begin{equation}\label{eq:qeehastings}
    h_Q(\mathbf B)=\min_{\substack{V<\mathbb C^n,\dim V\leq n/2}} \frac{\langle \Id-P_V,\Phi_{\mathbf B}(P_V)\rangle}{d\dim V}
\end{equation}
where $V<\mathbb C^n$ means that $V$ is a subspace of $\mathbb C^n$ %, $\dim V$ stands for the dimension of $V$  
and $P_V$ is the projector onto $V$ \cite{hastings2007random}. Note that if $\mathbf B$ is a permutation tuple determined by a graph $G$ and we restricted the minimum in \cref{eq:qeehastings} to subspaces spanned by vectors of the canonical basis we would recover the edge expansion of $G$. This is just one of the notions of quantum expansion introduced in the last two decades, which have found multiple applications in cryptography, learning theory and information theory \cite{hastings2007random,ben2007quantum,hastings2008classical,ambainis2009nonmalleable,kwok2019spectral,franks2020rigorous}. 

Another analogue of edge expansion in the quantum setting was recently introduced by Li, Qiao, Wigderson, Wigderson and Zhang \cite{li2022linear}. Given a bistochastic tuple $\mathbf B\in M_n(\mathbb C)^d$, its \emph{dimension edge expansion} is defined via 
\begin{equation}\label{eq:qeeLi}
    h_D(\mathbf B)=\min_{\substack{V<\mathbb C^n,\dim V\leq n/2}} \frac{\sum_{i\in [n]}\rank(B_i|_{V^\perp,V})}{d\dim V},
\end{equation}
where $B_i|_{V^\perp,V}=P_VB_i(\Id-P_V).$ Again, if $\mathbf B$ is a permutation tuple determined by a graph $G$ and we restricted the minimum in \cref{eq:qeeLi} to subspaces spanned by vectors of the canonical basis, we would recover the edge expansion of $G$. Li et al. showed that $d \cdot h_D\geq h_Q$ and they exhibited a sequence $\mathbf B_m$ of bistochastic $d$-tuples such that $\inf h_{D}(\mathbf B_m)>0$, but $\inf h_Q(\mathbf B_m)=0$ \cite[Theorems 1.8 and 1.9]{li2022linear}. To do that they considered a different expression of $h_Q(\mathbf B)$ in terms of the Schatten-2 norm\footnote{For $p\in [1,\infty)$ the Schatten-$p$ norm of a matrix is the $\ell_p$ norm of its singular values.}
\begin{equation*}
    h_Q(\mathbf B)=\underbrace{\min_{\substack{V<\mathbb C^n,\dim V\leq n/2}} \frac{\sum_{i\in [n]}\norm{B_i|_{V^\perp,V}}_{S_2}^2}{d\dim V}}_{:=h_{S_2}(\mathbf B)}.
\end{equation*}
If one substitutes $2$ in expression above by other $p\in [1,\infty)$  other natural  notions of expansion, namely 
\begin{equation*}
    h_{S_p}(\mathbf B)=\min_{\substack{V<\mathbb C^n,\dim V\leq n/2}} \frac{\sum_{i\in [n]}\norm{B_i|_{V^\perp,V}}_{S_p}^p}{d\dim V}.
\end{equation*}
These $S_p$ notions of expansion behave in a similar way as the case of $p=2$, as $d^{p/2} \cdot h_D\geq h_{S_p}$ and the aforementioned sequence $\mathbf B_m$ of bistochastic $d$-tuples also satisfies $\inf h_{D}(\mathbf B_m)>0$, but $\inf h_{S_p}(\mathbf B_m)=0$. This led Li et. al to ask whether all the $S_p$ notions of expansions are equivalent; $h_{S_p}$ and $h_{S_q}$ are equivalent if for every sequence $\mathbf B_m$ of bistochastic tuples $d$-tuples $\inf h_{S_p}(\mathbf B_m)>0$ if and only if $\inf h_{S_q}(\mathbf B_m)>0$. We give an affirmative answer to this question.

\begin{theorem}\label{theo:SpExpEquiv}
    Let $\mathbf{B}\in M_n(\mathbb C)^d$ be a bistochastic tuple. Let $p\geq q\geq 1$. Then, 
    \begin{enumerate}
        \item $h_{S_p}(\mathbf{B})\leq d^{\frac{p-q}{2}}h_{S_q}(\mathbf{B}),$\label{thm:itemi}
        \item $h_{S_p}(\mathbf{B})\geq [h_{S_q}(\mathbf{B})]^{\frac{p}{q}},$\label{thm:itemii}
        %\item $h_{S_p}(\mathbf{B})\geq \frac{1}{d}[h_{S_q}(\mathbf{B})]^{\frac{p}{q}},$\label{thm:itemii}
        %\item $[h_{S_p}(\mathbf B)]^{1/p}\geq h_{S_1}(\mathbf B).$\label{thm:itemiii}
    \end{enumerate}
    In particular, $h_{S_p}$ and $h_{S_q}$ are equivalent for every $p,q\in [1,\infty)$.
\end{theorem}

Li et al. pointed out that \cref{theo:SpExpEquiv} could be a quantum analogue of the equivalence between certain notions of graph expansion based on $\ell_p$ norms, a result by Matoušek \cite{matouvsek1997embedding}, who used this equivalence to prove that the metric spaces defined by graph expanders are hard to embed into $\ell_p$. However, we believe the analogy between our result and the one of Matoušek to be more notational than operational, so we suspect that \cref{theo:SpExpEquiv} will not help to construct finite metric spaces hard to embed into $S_p$. Nevertheless, we consider the latter an interesting open problem and propose a route towards solving it. We discuss this in \cref{sec:towardsSpembeddings}.

\section{All $S_p$ notions of expansion are equivalent}
\begin{proof}[ of~\Cref{theo:SpExpEquiv}.\ref{thm:itemi}]
    As $B_i$ belongs to a $d$-stochastic tuple, it satisfies  $B_i^*B_i=d\cdot\Id-\sum_{j\neq i} B_j^{*}B_j\preceq d\cdot\Id$, so all the singular values $(s_1^i,\dots,s_n^i)$ of $B_i|_{V^\perp,V}$ are at most~$\sqrt{d}$. Hence, 
	\begin{align*}
		\frac{\norm{B_i|_{V^\perp,V}}^p_{S_p}}{d^{\frac{p}{2}}}= \sum_{l\in [n]}\left(\frac{s_l^i}{\sqrt{d}}\right)^p\leq \sum_{l\in [n]}\left(\frac{s_l^i}{\sqrt{d}}\right)^q= \frac{\norm{B_i|_{V^\perp,V}}_{S_q}^q}{d^{\frac{q}{2}}},
	\end{align*}
	where in the inequality we have used that all $s_l^i/\sqrt{d}$ are at most 1 and that $p\geq q$. From here, the desired inequality follows immediately. 
\end{proof}
\begin{proof}[ of~\Cref{theo:SpExpEquiv}.\ref{thm:itemii}]
    \begin{comment}
    First we note that 
	\begin{align*}
		h_{S_p}(\mathbf{B})&\geq \frac{1}{d}\inf_{ \dim V\leq n/2} \max_{i\in [d]} \frac{\norm{B_i|_{V^\perp,V}}_{S_p}^p}{\dim V},\\
		h_{S_q}(\mathbf{B})&\leq\inf_{ \dim V\leq n/2} \max_{i\in [d]} \frac{\norm{B_i|_{V^\perp,V}}_{S_q}^q}{\dim V}.
	\end{align*}
	Thus it suffices to prove that for every matrix $B\in M_n(\mathbb C)$ and every $V<\mathbb{C}^n$ we have that 
	\begin{equation*}
		\frac{\norm{B|_{V^\perp,V}}_{S_p}^p}{\dim V}\geq 
		\left(\frac{\norm{B|_{V^\perp,V}}_{S_q}^q}{\dim V}\right)^{\frac{p}{q}}.
	\end{equation*}
	Indeed, given that $B|_{V^\perp,V}$ has rank at most $\dim V$, then it has at most $\dim V$ non-zero singular values $s=(s_1,\dots,s_{\dim V})$. Hence,  
	\begin{align*}
		\frac{\norm{B|_{V^\perp,V}}_{S_p}^p}{
		\dim V}=\norm{s}_{L_p^{\dim V}}^p\geq \norm{s}_{L_q^{\dim V}}^p=\left(\frac{\norm{B|_{V^\perp,V}}_{S_q}^q}{
		\dim V}\right)^{\frac{p}{q}},
	\end{align*} 
	where $L_p^{\dim V}$ stands for the $L_p$ space of dimension $\dim V$ with the uniform measure, and in the inequality we have used that $p\geq q$.
%\end{proof}
%\begin{proof}[ of~\Cref{theo:SpExpEquiv}.\ref{thm:itemiii}]
\end{comment}
    Let $(s^i_1,\ldots,s^i_{\dim V})$ be the singular values of $B_i|_{V^\perp,V}$.  Then by Hölder's inequality, 
	\begin{equation}\label{eq:firststepitemiii}
		\begin{split}
		\norm{B_i|_{V^\perp,V}}_{S_q}^q&=\sum_{l\in [\dim V]}|s_l^i|^q\leq\left(\sum_{l\in [\dim V]}|s_l^i|^p\right)^{q/p}\dim V^{1-\frac{q}{p}}=\norm{B_i|_{V^\perp,V}}_{S_p}^{q}\dim V^{1-\frac{q}{p}}.
		\end{split}
	\end{equation}
    We conclude that for every subspace $V<\mathbb{C}^n$, we have that
    \begin{equation*}
    	\begin{split}
    	\sum_{i\in [d]} \frac{\norm{B_i|_{V^\perp,V}}_{S_q}^q}{d\dim V}\leq 
    	\sum_{i\in [d]} \frac{\norm{B_i|_{V^\perp,V}}_{S_p}^{q}}{d\dim V^{q/p}}= d^{\frac{q}{p}-1}\sum_{i\in [d]}\left(\frac{\norm{B_i|_{V^\perp,V}}_{S_p}^p}{d\dim V}\right)^{q/p}
        \leq \left(\sum_{i\in [d]}\frac{\norm{B_i|_{V^\perp,V}}_{S_p}^p}{d\dim V}\right)^{q/p},
    	\end{split}
    \end{equation*}
	where in the first inequality we used \cref{eq:firststepitemiii} and in the second Hölder's inequality. From here, the desired inequality follows immediately.
\end{proof}	

\section{How much distortion is needed to embed metric spaces into $S_p$?}\label{sec:towardsSpembeddings}

\subsection{Embeddings into $\ell_p$ spaces}
A metric space $(M,\rho)$ can be embedded into $\ell_p$ with distortion $D$ if there is a map $f:M\to\ell_p$ such that
\begin{equation}\label{eq:distortion}
    \frac{1}{D}\rho(x,y)\leq\|f(x)-f(y)\|_{\ell_p}\leq\rho(x,y),\quad\forall x,y\in M.
\end{equation}
Finding the minimum distortion $D_{n,\ell_p}$ needed to embed any $n$-point metric space into $\ell_p$ was a long-standing question. The case $p=1$ was fully resolved by Bourgain (who proved the upper bound) and Linial, London and Rabinovich (who showed the lower bound). Together they concluded that $D_{\ell_1}=\Theta(\log n)$ \cite{bourgain1985lipschitz,linial1995geometry}.
The key ingredient of the lower bound of Linial et al. is that for a $d$-regular graph $G$ it holds that 
\begin{equation}\label{eq:h1>h}
    \underbrace{\inf_{f:G\to\ell_1}\frac{\frac{1}{|E|}\sum_{(i,j)\in E}\|f(i)-f(j)\|_{\ell_1}}{\frac{1}{n^2}\sum_{i,j\in [n]}\|f(i)-f(j)\|_{\ell_1}}}_{:=h_{\ell_1}(G)}\geq h(G).
\end{equation}
Matoušek realized that generalizing this result to other $\ell_p$ norms would allow one to prove similar lower bounds. With that purpose he (implicitly) introduced notions of graph expansion based on other $\ell_p$ norms, namely 
\begin{equation*}
    h_{\ell_p}(G):=\min_{f:G\to\ell_p}\left(\frac{\frac{1}{|E|}\sum_{(i,j)\in E}\|f(i)-f(j)\|_{\ell_p}^p}{\frac{1}{n^2}\sum_{i,j\in [n]}\|f(i)-f(j)\|_{\ell_p}^p}\right)^{1/p}.
\end{equation*}
Then, he showed that 
\begin{equation}\label{eq:hlp>hl1}
    h_{\ell_p}(G)\geq \frac{h_{\ell_1}(G)}{4dp},
\end{equation}
which, together with \eqref{eq:h1>h} and the existence of a sequence $G_n$ of $n$-vertex $d$-regular graphs such that $\inf h(G_n)=C>0$, implies that 
\begin{equation}\label{eq:hlp>0}
    \inf h_{\ell_p}(G_n)\geq \frac{C}{4dp}.
\end{equation} 
From \cref{eq:hlp>0}, Matoušek showed with an elegant argument that the metric spaces defined by these graphs require distortion $D_{n}=\Omega(\log (n)/p)$ to be embedded into $\ell_p$. Indeed, consider the metric spaces $([n],\rho_n)$, where $
\rho_n(i,j)$ is the length of the shortest path in $G_n$ from $i$ to $j$. Define $$R_{\rho_n
}=\left(\frac{\frac{1}{|E|}\sum_{(i,j)\in E}\rho_n(i,j)^p}{\frac{1}{n^2}\sum_{i,j\in [n]}\rho_n(i,j)^p}\right)^{1/p}.$$ The numerator of $R_{\rho
_n}$ is 1, while the denominator is $\Omega(\log_d n
)$, because there are at most $d^1+\dots+d^{\log_d(n/4)}\leq 2d^{\log_d(n/4)}=n/2$ pairs of vertices whose shortest-path distance is at most $\log_d(n/4).$ Hence, treating $d$ as a constant, it follows that  $R_{\rho_n}=O(1/\log (n))$. Now let $f_n$ be embeddings of $([n],\rho_n)$ into $\ell_p$ with distortion $D_n$ and define $R_{f_n}$ in an analogue way. As $f_n$ are embeddings with distortion $D_n$ it follows that $R_{f_n}=O( D_n/\log (n))$ and \cref{eq:hlp>0} implies that $R_{f_n}= \Omega(1/p)$. Putting both bounds together yields $D_{n,\ell_p}=\Omega(\log (n)/p)$. Matoušek also proved the matching upper bound by mimicking the argument of Bourgain, giving $D_{n,\ell_p}=\Theta(\log (n)/p)$. 

\subsection{What about $S_p$?}
As pointed out by Li et al. \cite{li2022linear}, our \cref{theo:SpExpEquiv} can be regarded as a quantum analogue of \cref{eq:hlp>hl1} proved by Matoušek, in the sense that both prove relations between notions of expansion based on similar norms. However, we believe that this resemblance is notational, rather than operational, and so our \cref{theo:SpExpEquiv} does not seem to help to construct metric spaces that are hard to embed into $S_p$. To the best of our knowledge there is no literature about this problem, so we pose it as an open question. 
\begin{question}
    Let $p\in[1,\infty)$. What is the minimum distortion $D_{n,S_p}$ needed, as a function of $n$ and $p$, to embed any $n$-point metric space into $S_p$?
\end{question}
To answer this question it would suffice to show that the following family of notions of graph expansion based on $S_p$ norms for $p\in [1,\infty)$ are related in a certain way $$\tilde h_{S_p}(G):=\min_{f:G\to S_p}\left(\frac{\frac{1}{|E|}\sum_{(i,j)\in E}\|f(i)-f(j)\|^p_{S_p}}{\frac{1}{n^2}\sum_{i,j\in [n]}\|f(i)-f(j)\|^p_{S_p}}\right)^{1/p}.$$
Indeed, as $S_2$ is isometric to $\ell_2$, we have that the there exists a sequence of $n$-vertex $d$-regular graphs $G_n$ such that $\inf \tilde h_{S_2}(G_n)>0$. Hence, if one manages to show that for every $d$-regular graph $G$ it holds that $$\tilde h_{S_p}(G)=\Omega_d\left(\frac{\tilde h_{S_2}(G)}{p}\right),$$ we would have that $\inf \tilde h_{S_p}(G_n)=\Omega(1/p)$, from which we could reproduce the argument of Matoušek and prove that $D_{n,S_p}=\Omega(\log(n)/p)$. The matching upper bound follows from the fact that $\ell_p$ can be embedded isometrically into $S_p$.

\paragraph{Acknowledgements.}We thank Jop Bri\"et for detailed comments and discussions. We thank Davi Castro-Silva, Yinan Li, Youming Qiao,  Yuval Wigderson and Chuanqi Zhang for useful comments and discussions.  This research was supported by the European Union’s Horizon 2020 research and innovation programme under the Marie Sk{\l}odowska-Curie grant agreement no. 945045, and by the NWO Gravitation project NETWORKS under grant no. 024.002.003. 
\bibliographystyle{alphaurl}
\bibliography{Bibliography}
\end{document}